\def\vers{Nov.~19, 2009, v.1}
\magnification=1200
\hsize=6.5truein
\vsize=8.9truein
\font\bigfont=cmr10 at 14pt
\font\mfont=cmr9
\font\sfont=cmr8
\font\mbfont=cmbx9

\def\nin{\noindent}
\def\bs{\bigskip}
\def\ms{\medskip}

\def\scirc{\,\raise.2ex\h{${\scriptstyle\circ}$}\,}
\def\simto{\buildrel\sim\over\longrightarrow}

\def\mtim{\h{$\times$}}
\def\mcup{\h{$\bigcup$}}
\def\mcopr{\h{$\coprod$}}
\def\a{{\alpha}}
\def\A{{\cal A}}
\def\b{\backslash}
\def\C{{\bf C}}

\def\DD{\Delta}
\def\G{{\Gamma}}
\def\GG{\widehat{\Gamma}}
\def\h{\hbox}
\def\H{{\bf H}}

\def\LL{\widehat{\cal L}}

\def\tnu{\widetilde{\nu}}
\def\O{{\cal O}}
\def\q{\quad}

\def\O{{\cal O}}
\def\Q{{\bf Q}}
\def\R{{\bf R}}
\def\V{{\cal V}}
\def\VV{\widehat{\cal V}}
\def\Z{{\bf Z}}
\def\Gr{\hbox{\rm Gr}}
\def\NF{\hbox{\rm NF}}
\def\Im{\hbox{\rm Im}}
\def\Ker{\hbox{\rm Ker}}

\def\Ext{\hbox{\rm Ext}}
\def\Hom{\hbox{\rm Hom}}
\def\MHS{{\rm MHS}}
\def\VMHS{{\rm VMHS}}
\def\BPS{{\rm BPS}}

\def\ms{\medskip}
\def\da{\downarrow}

\def\simto{\buildrel\sim\over\to}
\def\1{\hskip1pt}
\def\({{\rm (}}
\def\){{\rm )}}

\hbox{}
\vskip 1cm

\centerline{\bigfont Generalization of N\'eron models of Green,
Griffiths and Kerr}

\medskip
\centerline{(Joint work with P.~Brosnan and G.~Pearlstein)}

\bs
\centerline{\it Dedicated to Professor Sampei Usui}

\bs
\centerline{Morihiko Saito}

\ms
\centerline{RIMS Kyoto University, Kyoto 606-8502 Japan}

\bs\medskip
{\narrower\noindent
{\mbfont Abstract.} {\mfont
We explain some recent developments in the theory of N\'eron models
for families of Jacobians associated to variations of Hodge
structures of weight $-1$.}
\par}

\bs\bs
\centerline{\bf 1. Classical N\'eron models}

\bs\nin
{\bf 1.1.}
Let $\A$  be an abelian scheme over a smooth curve $S^*\subset S$.
There is a unique group scheme $\A_S$ over $S$, called the
N\'eron model, and satisfying the following property:
For any smooth $T$ over $S$, we have
$$\A_S(T)=\A(T_{S^*}),\q\h{i.e.}\q\A_T(T)=\A_{T_{S^*}}(T_{S^*}).$$
Let $\H$ be the variation of Hodge structure of level 1 and weight
$-1$ corresponding to $\A$. Then we have for $s\in S^*$
$$\A_s=J(H_s)\,\,(:=H_{s,\Z}\b H_{s,\C}/F^0H_{s,\C}).$$
Note that the right-hand side is isomorphic to
$$\Ext^1_{\MHS}(\Z,H_s),$$
i.e.\ its element corresponds to the short exact sequence of MHS
(see [Ca])
$$0\to H_s\to H'_s\to\Z\to 0,$$
where MHS denotes the abelian category of mixed $\Z$-Hodge
structures [D2].

Assume the monodromy is {\it unipotent} at $0\in S\setminus S^*$.
By [Sd] we have the limit mixed Hodge structure
$$H_{\infty}=\bigl((H_{\infty,\C};F,W),\,(H_{\infty,\Q},W),
\,H_{\infty,\Z}\bigr).$$
This is closely related to the N\'eron model.
Indeed, there is a short exact sequence
$$0\to\A_{S,0}^0\to\A_{S,0}\to G_0\to 0,$$
where
$$\A_{S,0}^0:=H_{\infty,\Z}^{\rm inv}\b H_{\infty,\C}/
F^0H_{\infty,\C},$$
and
$$\eqalign{G_0:=\,\,&H^1(\DD^*,\H_{\Z})_{\rm tor}
={\rm Coker}(T_{\Z}-id)_{\rm tor}\cr
=\,\,&(\Im(T_{\Q}-id)\cap H_{\infty,\Z})/\Im(T_{\Z}-id)\cr
=\,\,&\Ker(T_{\Q/\Z}-id)/({\rm Im\,\,of}\,\,\Ker(T_{\Q}-id)).}$$
For the last isomorphism we use the snake lemma applied to the
endomorphism $T$ of the short exact sequence
$$0\to\H_{\Z}\to\H_{\Q}\to\H_{\Q/\Z}\to 0.$$
This can be used to get a torsion normal function
corresponding to a torsion cohomology class.

\ms\nin
{\bf 1.2.~Example.} Let $\A$ be a family of elliptic curves
with monodromy $T=\pmatrix{1&r\cr 0&1}$.
Then
$$\A_{S,0}^0=\Z\b \C^2/\C=\C^*,\q
G_0=\Z/r\Z,\q \A_{S,0}=\mcopr^r\C^*.$$

\bs\bs
\centerline{\bf 2. Generalization by Zucker and Clemens}

\bs\nin
{\bf 2.1.~Generalization by Zucker.}
Let $\H$ be a variation of Hodge structure of weight $-1$ on $S^*$.
we have the family of Jacobians
$$J(\H):=\mcopr_{s\in S^*}J(H_s).$$
Let $\LL$ be the Deligne extension of $\H_{\O}$ over $S$ (see [D1]),
$\VV$ the vector bundle corresponding to $\LL/F^0\LL$,
and $\GG$ the image of $j_*\H_{\Z}$ in $\VV$ where
$j:S^*\hookrightarrow S$.

\ms\nin
{\bf 2.2.~Definition} (Zucker extension) [Zu].
$$J_S^Z(\H):=\GG\b \VV\,\,\h{(fiberwise).}$$

\ms\nin
Assume $\H$ geometric, i.e.
$\H=R^{2p-1}f_*\Z_{X^*}(p)$ with $f:X^*\to S^*$.
Then we have the following.

\ms\nin
{\bf 2.3.~Theorem} (El Zein, Zucker) [EZ]. Let $\nu$ be a normal
function defined by an algebraic cycle $\sigma$ with
$\gamma(\sigma|_{X_s})=0$.
Then $\nu$ extends to a section of $J_S^Z(\H)$ over $S,$ if
$\gamma(\sigma)=0.$

\ms
Here $\gamma(\sigma)$ denotes the cohomology class
as a cycle.

\ms\nin
{\bf 2.4.~Generalization by Clemens.} Assume

\ms\nin
{\it Hypothesis $(C)$} :
$N^2=0$ and $\Gr_0^WH_{\infty}$ has type (0,0).

\ms\nin
Then we have the following.

\ms\nin
{\bf 2.5.~Theorem}~(Clemens) [Cl].
{\it There is $J_S^C(\H)$ \(Clemens N\'eron model\1\) such that
any normal function $\nu$ on $S^*$ defined by an algebraic
cycle is extended to a section of $J_S^C(\H)$ over $S$.
Moreover, there is a short exact sequence
$$0\to J_S^Z(\H)_0\to J_S^C(\H)_0\to G_0\to 0,$$
with}
$$J_S^Z(\H)_0=H_{\infty,\Z}^{\rm inv}\b H_{\infty,\C}/
F^0H_{\infty,\C},$$
$$G_0=H^1(\DD^*,\H_{\Z})_{\rm tor}.$$

\ms
In fact, $J_S^C(\H)$ is obtained by gluing $J_S^Z(\H)$.

\bs\bs
\centerline{\bf 3. Improvement using admissible normal functions}

\bs\nin
{\bf 3.1.} By [Ca] we have a canonical isomorphism
$$\Ext_{\MHS}^1(\Z,\H_s)=J(\H_s)\,
(:=H_{s,\Z}\b H_{s,\C}/F^0H_{s,\C}),$$
and a normal function $\nu\in\NF(S^*,\H)$ (which is a holomorphic
section of $J(\H)$) corresponds to a short exact sequence
$$0\to\H\to\H'\to\Z_{S^*}\to 0,$$
where the Griffiths transversality of $\nu$ corresponds to that
of $\H'$.
So we get
$$\NF(S^*,\H)=\Ext_{\VMHS}^1(\Z_{S^*},\H)\,\,(\subset J(\H)(S^*)).$$
By construction [Ca], this is induced by taking the difference of
the two local splittings $\sigma_F$ and $\sigma_{\Z}$ of the above
short exact sequence, where the $\sigma_F$ is compatible with $F$
and $\sigma_{\Z}$ is defined over $\Z$.
Note that the cohomology class $\gamma(\nu)\in H^1(S^*,\H)$ is
defined by using the cohomology long exact sequence of the above
short exact sequence.

Let $S^*\subset S$ be a partial compactification such that
$S\setminus S^*\subset S$ is closed analytic.

\ms\nin
{\bf 3.2.~Definition} (Admissible normal functions with respect to
$S^*\subset S$, see [Sa2]).
$$\NF(S^*,\H)^{\rm ad}_S:=\Ext_{\VMHS(S^*)^{\rm ad}_S}^1
(\Z_{S^*},\H)\,\,(\subset J(\H)(S^*)),$$
where $\VMHS(S^*)^{\rm ad}_S$ denotes the category of admissible
variation of mixed Hodge structure with respect to $S^*\subset S$.

\ms\nin
{\bf 3.3.~Definition.} The category of admissible VMHS in the
one-dimensional case is defined by the following two conditions of
Steenbrink and Zucker [SZ] where $(S,S^*)=(\DD,\DD^*)$:

\ms\nin
(a) The $\Gr^p_F\Gr^W_k\LL$ are free $\O_{\DD}$-modules
(in the unipotent case).

\ms\nin
(b) The relative monodromy filtration exists.

\ms\nin
{\bf 3.4.~Remarks.} (i) In the non-unipotent case, (a) is not
sufficient, and we have to take a ramified covering to reduce
to the unipotent monodromy case (or use the $V$-filtration).

\ms\nin
(ii) Condition (b) in the weight $-1$ case is equivalent to a
splitting of the short exact sequence of $\Q$-local systems over
$\DD^*$.

\ms\nin
(iii) The generalization to the higher dimensional
case is by the curve test, see [Ka].

\ms\nin
{\bf 3.5.~Remark.} Taking a multivalued lifting $\tnu$ of
$\nu$ in $\V$,
conditions (a), (b) correspond to the conditions given
by Green, Griffiths, Kerr [GGK] (in the unipotent case):

\ms\nin
(a)$'$ $\tnu$ has a logarithmic growth.

\ms\nin
(b)$'$ $T\1\tnu-\tnu\in\Im(T_{\Q}-id)$.

\ms\nin
Note that the variation $T\1\tnu-\tnu$ gives the cohomology class
of $\nu$.

\ms
Using the theory of admissible normal functions, we can show the
following.

\ms\nin
{\bf 3.6.~Theorem} [Sa2]. {\it Theorem~{\rm (2.5)} holds for any
admissible normal function $\nu$ \(not necessarily associated to
an algebraic cycle\) without assuming the hypothesis $(C)$,
and $J_S^C(\H)$ has a structure of a complex Lie group over $S$.}

\ms
The key point is the following generalization of Theorem~(2.3).

\ms\nin
{\bf 3.7.~Proposition} [Sa2]. {\it For an admissible normal function
$\nu$ on $S^*$, $\nu$ extends to a section of $J_S^Z(\H)$ if and
only if $\gamma_0(\nu)=0$.}

\ms
Here $\gamma_0(\nu)$ is the cohomology class of $\nu|_{\DD^*}$
where $\DD\subset S$ is a sufficiently small disk with center 0.
Theorem~(3.6) is then proved by setting
$$J_S^C(\H)|_{\DD}:=\mcup_{g\in G_0}\,
\bigl(\nu_g+J_S^Z(\H)|_{\DD}\bigr),$$
where $\nu_g\in\NF(\DD^*,\H_{\DD^*})_{\DD}^{\rm ad}$
such that $\gamma_0(\nu_g)=g$.
This is independent of the choice of $\nu_g$ by Proposition~(3.7).

\ms\nin
{\bf 3.8.~Remark.}
If $\H$ corresponds to an abelian scheme $\A$, then
$$\A_S^{\rm an}\simto J_S^C(\H),$$
even in the non-unipotent case, see [Sa2], 4.5.

\bs\bs
\centerline{\bf 4. Generalization by Green, Griffiths and Kerr
{\rm [GGK]}}

\bs\nin
{\bf 4.1.~Problem.}
In general, $J_S^Z(\H),\,J_S^C(\H)$ are not Hausdorff as is shown
by an example in [Sa2], 3.5(iv) where $H_{\infty}$ has type
$$(1,-1),(-1,1),(0,-2),(-2,0).$$

\ms\nin
{\bf 4.2.~Theorem} (Green, Griffith, Kerr) [GGK].
{\it Except for the last assertion on the structure of a complex
Lie group over $S$, Theorem~{\rm (3.6)} also holds for a subspace
$J_S^{\rm GGK}(\H)$ of $J_S^C(\H)$ which is obtained by replacing
$J_S^Z(\H)_0$ with
$$J_S^{\rm GGK}(\H)_0^0:=J(H_{\infty}^{\rm inv})\,
(=H_{\infty,\Z}^{\rm inv}\b H_{\infty,\C}^{\rm inv}/
F^0H_{\infty,\C}^{\rm inv}),$$
where $H_{\infty}^{\rm inv}:=\Ker\,N\subset H_{\infty}$ with
$N=\log T$.}

\ms
Here the monodromy is assumed {\it unipotent}.
Note that $J_S^{\rm GGK}(\H)$ is not an analytic space in the
usual sense, see [GGK].
We have moreover the following.

\ms\nin
{\bf 4.3.~Theorem} [Sa3]. {\it As a topological space endowed with
the quotient topology, $J_S^{\rm GGK}(\H)$ is Hausdorff
\(assuming $\dim S=1$\).}

\ms\nin
{\bf 4.4~Remark.}
These assertions can be extended to the case
$\dim S>1$ if $D:=S\setminus S^*$ is {\it smooth}.

\ms\nin
{\bf 4.5.~Corollary.} {\it The closure of the zero locus of an
admissible norma function is analytic if $D$ is smooth.}

\ms\nin
{\bf 4.6.~Remark.} This is independently proved by P.~Brosnan
and G.~Pearlstein using another method, and they recently
give a proof in the general case [BP], see also [KNU], [Sl].
The generalization of (4.3) seems to be closely related to [CKS].

\bs\bs
\centerline{\bf 5. Generalization by Brosnan, Pearlstein and Saito
{\rm [BPS]}}

\bs\nin
{\bf 5.1.}~Assume $S^*$ smooth, but $S$ may be singular.
Consider the inclusions $j:S^*\hookrightarrow S$ and
$i_s:\{s\}\hookrightarrow S$ for $s\in S$. Set
$$\eqalign{H_s&:=H^0i_s^*(\R j_*\H)\,
(=\lim_{U\ni 0} H^0(U\cap S^*,\H)),\cr
J(H_s)&:=\Ext^1(\Z,H_s)\,(=H_{s,\Z}\b H_{s,\C}/F^0H_{s,\C}).}$$

\ms\nin
{\bf 5.2.~Identity component of the BPS-N\'eron model} [BPS].
Define
$$J^{\BPS}_S(\H)^0:=\mcopr_{s\in S}J(H_s)\,\,
\h{(set-theoretically)}.$$
Then $J^{\BPS}_S(\H)^0$ is a topological space
(using a resolution of $S$), and
$$J^{\BPS}_S(\H)^0_{S_{\alpha}}:=J^{\BPS}_S(\H)^0|_{S_{\alpha}}$$
is a Lie group over $S_{\alpha}$ for any stratum $S_{\alpha}$ of a
Whitney stratification of $(S,D)$.
Moreover, $\nu$ defines a continuous section of $J^{\BPS}_S(\H)^0$
if $\gamma_s(\nu)=0\,\,(\forall\,s\in D)$.
Indeed, $\nu$ corresponds to a short exact sequence
$$0\to\H\to\H'\to\Z_{S^*}\to 0,$$
and this induces a long exact sequence
$$0\to H^0i_s^*\R j_*\H\to H^0i_s^*\R j_*\H'\to\Z\to H^1i_s^*
\R j_*\H.$$
So we have $\nu_s\in J(H_s)$ if $\gamma_s(\nu)=0$.

For the reduction to the normal crossing case, set
$${\rm NF}_{(s)}(S^*,\H)_S^{\rm ad}:=
\{\nu\in{\rm NF}(S^*,\H)_S^{\rm ad}\mid\gamma_s(\nu)=0\}.$$
Then we have the following.

\ms\nin
{\bf 5.3.~Proposition.} {\it For $\pi:S'\to S$ proper with
$S'{}^*:=\pi^{-1}(S^*)$ smooth, we have the commutative diagram}
$$\matrix{{\rm NF}_{(s)}(S^*,\H)_S^{\rm ad}&\to&
{\rm NF}_{(s')}(S'{}^*,\pi^*\H)_{S'}^{\rm ad}\cr
\da&\raise12pt\hbox{ }\raise-8pt\hbox{ }&\da\cr
J(H_s)&\to&J(H_{s'})}$$

\ms\nin
{\bf 5.4.~Definition.} Set $G_s:=\{\h{images of admissible
$\nu$}\}\subset H^1i_s^*\R j_*\H$, and
$$G:=\mcopr_{s\in S}G_s.$$
Then we can prove the following.

\ms\nin
{\bf 5.5.~Theorem} [BPS]. {\it There exists $J^{\BPS}_S(\H)$ over
$S$ together with an exact sequence
$$0\to J^{\BPS}_S(\H)^0\to J^{\BPS}_S(\H)\to G\to 0,$$
and any admissible $\nu$ on $S^*$ is extendable to a
continuous section of $J^{\BPS}_S(\H)$ over $S$.}

\ms
This is shown by generalizing the gluing argument in the
proof of Theorem~(3.6) in the one-dimensional case.

\ms\nin
{\bf 5.6.~Remarks.} (i) A.~Young [Yo] constructed a generalization
of N\'eron model for families of Abelian varieties defined
on the complement of a divisor with normal crossings
where he assumes that the local monodromies are unipotent and
the identity component is similar to the classical construction [Na].

\ms\nin
(ii)
C.~Schnell [Sl] has given a definition of a N\'eron model whose
identity component $J_S^{\rm Sch}(\H)^0$ is Hausdorff by using
the Hodge filtration of Hodge modules [Sa1] where the partial
compactification $S$ is smooth although $S\setminus S^*$ is not
necessarily a divisor with normal crossings (see [SS] for the
one-dimensional case).

\ms\nin
(iii) In the case of abelian schemes over curves,
we get something different from the classical N\'eron
model if the monodromy is {\it non-unipotent}.
Indeed, for a family of elliptic curves with non-unipotent monodromy,
we have
$$J(H_0)=pt,$$ (since $H_0:=H^0i_0^*\H=0$), and
there is a `{\it blow-down}' map
$$\A_S^{\rm an}= J_S^C(\H)\to J^{\BPS}_S(\H),$$
such that the image of $\A_{S,0}^{\rm an}=J_S^C(\H)_0=\C$ is
$J(H_0)=pt$.

\ms\nin
(iv) It is very difficult to determine $G_0$
even in the normal crossing case.
Indeed, there is a cohomology class map
$$\NF((\DD^*)^n,\H_{\Q})^{\rm ad}_{\DD^n}\to
\Hom_{\MHS}(\Q,{\cal H}^1({\rm IC}_{(\DD^*)^n}\H_{\Q})_0),$$
which is surjective if $H$ is a nilpotent orbit,
but it can be non-surjective in general [Sa4].
However, this does not seem to contradicts the strategy of
Green, Griffith for solving the Hodge conjecture
since it seems to occur only in rather artificial occasions,
e.g.\ when an unnecessary blowing-up is made.

\ms\nin
(v) As is remarked by C.~Schnell [Sl], the topology of the N\'eron
models which graph any admissible normal functions with non-vanishing
cohomology classes can be rather complicated even in the abelian
scheme case [Yo] as is shown by the example below if $\dim S\ge 2$.

\ms\nin
{\bf 5.7.~Example.}
Let $S=\DD^2$, $S^*=(\DD^*)^2$, and $\rho:S\to\DD$ be the morphism
defined by $\rho(t_1,t_2)=t_1t_2$.
Let $\H^1$ be the nilpotent orbit of weight $-1$ having an integral
basis $e_0,e_1$ such that $Ne_0=e_1$ and $F^0H_{\infty,\C}=\C\1 e_0$.
Set $\H=\rho^*\H^1$. Then
$$\eqalign{&G_0\cong\Z\,\,\h{(non-canonically)},\q
G_s=0\,\,(s\ne 0),\cr
&J_S^{\BPS}(\H)^0=\C^*\mtim\DD^2/\G'\q\h{with}\cr
&\G'=\{(x,t_1,t_2)\in\C^*\mtim(\DD^*)^2\mid
x=(t_1t_2)^k\,(k\in\Z)\}.}$$
For $(p,\a)\in\Z\mtim\C^*$, there is an admissible normal function
$\nu_{p,\a}$ with respect to $S^*\subset S$ defined by
$$x=\a t_1^p=\a t_1^{k+p}t_2^k\,\,\h{mod}\,\,\Gamma'
\,\,\h{for any}\,\,k\in\Z.$$
Its cohomology class $\gamma_0(\nu_{p,\a})\in G_0\cong\Z$ is equal
to $p$ up to a sign, and the closure of its graph over $S^*$
contains the $0$-component $J_S^{\BPS}(\H)^0_0=\C^*$ of
$J_S^{\BPS}(\H)_0$ if $|p|\ge 2$ (restricting over a curve defined
by $t_2=\beta\1t_1^{|p|-1}$ with $\beta\in\C^*$).
However, the extended section $\overline{\nu}_{p,\a}$ over $S$
passes through the $\gamma_0(\nu_{p,\a})$-component of
$J_S^{\BPS}(\H)_0$ by construction.
(Here the value $\overline{\nu}_{p,\a}(0)$ at the origin depends
also on $\a$, and this induces an isomorphism
$\Z\mtim\C^*\simeq J_S^{\BPS}(\H)_0$.
In particular, for any admissible normal function $\nu$, there is
a unique $(p,\a)\in\Z\mtim\C^*$ such that
$\overline{\nu}(0)=\overline{\nu}_{p,\a}(0)$.
Then the above assertion on the closure of the graph holds for any
admissible function $\nu$ with $|\gamma_0(\nu)|\ge 2$.)

The above argument implies that $J_S^{\BPS}(\H)$ cannot be
Hausdorff as a topological space, and this would be the same for
the N\'eron models in [Sl], [Yo] which should coincide with
$J_S^{\BPS}(\H)$ for this simple example (even if the identity
component $J_S^{\rm Sch}(\H)^0$ in [Sl] is Hausdorff in general).

\vfill\eject
\centerline{{\bf References}}

\medskip
{\sfont
\item{[BP]}
P.~Brosnan and G.~Pearlstein,
On the algebraicity of the zero locus of an admissible normal
function, preprint (arXiv:0910.0628)

\item{[BPS]}
P.~Brosnan, G.~Pearlstein and M.~Saito,
A generalization of the N\'eron models of Green, Griffiths and Kerr,
preprint (arXiv:0809.5185)

\item{[Ca]}
J.~Carlson, Extensions of mixed Hodge structures, in Journ\'ees de
G\'eom\'etrie Alg\'ebrique d'Angers 1979, Sijthoff-Noordhoff Alphen
a/d Rijn, 1980, pp.~107--127.

\item{[CKS]}
E.~Cattani, A.~Kaplan and W.~Schmid,
Degeneration of Hodge structures, Ann.\ of Math.\ 123 (1986),
457--535. 

\item{[Cl]}
H.~Clemens, The N\'eron model for families of intermediate Jacobians
acquiring ``algebraic'' singularities, Publ.\ Math.\ IHES 58 (1983),
5--18.

\item{[D1]}
P.~Deligne, Equations diff\'erentielles\`a points singuliers
r\'eguliers, Lect. Notes in Math. vol. 163, Springer, Berlin, 1970.

\item{[D2]}
P.~Deligne, Th\'eorie de Hodge II, Publ. Math. IHES 40 (1971),
5--58.

\item{[EZ]}
F.~El Zein and S.~Zucker, Extendability of normal functions associated
to algebraic cycles, in Topics in transcendental algebraic geometry,
Ann. Math. Stud., 106, Princeton Univ. Press, Princeton, N.J., 1984,
pp.~269--288.

\item{[GGK]}
M.~Green, P.~Griffiths and M.~Kerr,
N\'eron models and limits of Abel-Jacobi mappings (preprint).

\item{[Ka]}
M.~Kashiwara, A study of variation of mixed Hodge structure,
Publ.\ RIMS, Kyoto Univ. 22 (1986), 991--1024.

\item{[KNU]}
K.~Kato, C.~Nakayama and S.~Usui,
Moduli of log mixed Hodge structures (arXiv:0910.4454).

\item{[Na]}
Y.~Namikawa,
A new compactification of the Siegel space and degeneration of
Abelian varieties, I and II, Math.\ Ann.\ 221 (1976), 97--141 and
201--241.

\item{[Sa1]}
M.~Saito, Mixed Hodge modules, Publ.\ RIMS, Kyoto Univ. 26
(1990), 221--333.

\item{[Sa2]}
M.~Saito,
Admissible normal functions, J.\ Algebraic Geom.\ 5 (1996),
235--276.

\item{[Sa3]}
M.~Saito,
Hausdorff property of the N\'eron models of Green, Griffiths and Kerr
(arXiv:0803.2771).

\item{[Sa4]}
M.~Saito,
Cohomology classes of admissible normal functions, preprint
(arXiv:0904.1593).

\item{[SS]}
M.~Saito and C.~Schnell,
A variant of N\'eron models over curves, preprint (arXiv:0909.4276).

\item{[Sd]}
W.~Schmid, Variation of Hodge structure: The singularities of the
period mapping, Inv.\ Math.\ 22 (1973), 211--319.

\item{[Sl]}
C.~Schnell, Complex analytic N\'eron models for arbitrary families
of intermediate Jacobians, preprint (arXiv:0910.0662).

\item{[SZ]}
J.H.M.~Steenbrink and S.~Zucker, Variation of mixed Hodge structure,
I, Inv.\ Math.\ 80 (1985), 489--542.

\item{[Yo]}
A.~Young, Complex analytic N\'eron models for degenerating Abelian
varieties over higher dimensional parameter spaces,
Ph.\ D.\ Thesis, Princeton University, 2008.

\item{[Zu]}
S.~Zucker, Generalized intermediate Jacobians and the theorem on
normal functions, Inv.\ Math.\ 33 (1976),185--222.

\ms\nin
\vers
\bye